# Asymptotics of the Euler transform of Fibonacci numbers

Václav Kotěšovec

e-mail: kotesovec2@gmail.com

August 7, 2015

**Abstract:** The generating function for the sequence A166861 in the OEIS ("Euler transform of Fibonacci numbers") is

$$\prod_{k=1}^{\infty} \frac{1}{(1-x^k)^{F_k}}$$

where $F(k)$ are the Fibonacci numbers (A000045). This paper analyzes the more general generating function

$$U(x) = \prod_{k=1}^{\infty} \frac{1}{(1-x^k)^{F_{k+z}}}$$

where $z$ is a nonnegative integer, which provides asymptotics for the sequences A166861 (z=0), A200544 (z=1) and A260787 (z=2) in the OEIS.

**Main result:**

$$a_n \sim \frac{\varphi^{n+\frac{z}{4}} \, e^{\left(\frac{\varphi}{10} - \frac{1}{2}\right) F_z - \frac{1}{10} F_{z+1} + \frac{2\,\varphi^{z/2}}{5^{1/4}} \sqrt{n} + S}}{2\sqrt{\pi} \, 5^{1/8} \, n^{3/4}}$$

where

$$S = \sum_{k=2}^{\infty} \frac{F_z + F_{z+1}\,\varphi^k}{(\varphi^{2k} - \varphi^k - 1)\,k}$$

and

$$\varphi = \frac{1+\sqrt{5}}{2}$$

is the golden ratio (A001622).

**Proof**:

We have the Maclaurin series

$$\log(1-x) = -\sum_{k=1}^{\infty} \frac{x^k}{k}$$

and so

$$\log(U(x)) = \log\left(\prod_{j=1}^{\infty} \frac{1}{(1-x^j)^{F_{j+z}}}\right) = -\sum_{j=1}^{\infty} F_{j+z} \log(1-x^j) = \sum_{j=1}^{\infty} F_{j+z} \sum_{k=1}^{\infty} \frac{x^{jk}}{k} = \sum_{k=1}^{\infty} \frac{x^k(F_z x^k + F_{z+1})}{k(1-x^k-x^{2k})}$$

Using the saddle-point method, see [2], equation (12.9), we have

$$a_n \sim \frac{U(r_n)}{\sqrt{2\pi * b(r_n)} * r_n^n}$$

The saddle-point equation is

$$r_n * \frac{U'(r_n)}{U(r_n)} = n$$

$$x * \frac{U'(x)}{U(x)} = x * \frac{d}{dx}\log(U(x)) = x * \frac{d}{dx}\sum_{k=1}^{\infty} \frac{x^k(F_z x^k + F_{z+1})}{k(1-x^k-x^{2k})} = \sum_{k=1}^{\infty} \frac{x^k(F_{z+1}(x^{2k}+1) - F_z x^k(x^k-2))}{(x^{2k}+x^k-1)^2}$$

$$\sum_{k=1}^{\infty} \frac{r_n^k(F_{z+1}(r_n^{2k}+1) - F_z r_n^k(r_n^k-2))}{(r_n^{2k}+r_n^k-1)^2} = n$$

For an asymptotic solution set $k = 1$ and the dominant root is then

$$r_n = \varphi - 1 - \frac{\varphi^{\frac{z}{2}-1}}{5^{1/4}\sqrt{n}} + \frac{\varphi^{z-1}}{2\sqrt{5}\,n} + O\left(\frac{1}{n^{3/2}}\right)$$

where $\varphi = \frac{1+\sqrt{5}}{2}$ is the golden ratio. It is important to note that taking only two terms the asymptotic expansion $(\varphi - 1) - \frac{\varphi^{z/2-1}}{5^{1/4}\sqrt{n}}$ is insufficient, three terms are needed. An eventual term $n^{-3/2}$ can be ignored.

Now we compute

$$\frac{1}{r_n^n} \sim \frac{1}{\left(\varphi - 1 - \frac{\varphi^{\frac{z}{2}-1}}{5^{1/4}\sqrt{n}} + \frac{\varphi^{z-1}}{2\sqrt{5}\,n}\right)^n} \sim \varphi^n \exp\left(\frac{\varphi^{z/2}\sqrt{n}}{5^{1/4}}\right)$$

$$b(x) = \frac{x\,U'(x)}{U(x)} + \frac{x^2\,U''(x)}{U(x)} - \frac{x^2\,U'(x)^2}{U(x)^2} = \frac{x\,U'(x)}{U(x)} + x^2\left(\frac{d}{dx}\right)^2 \log(U(x))$$

We obtain

$$b(x) = \frac{x\,U'(x)}{U(x)} + x^2\left(\frac{d}{dx}\right)^2 \sum_{k=1}^{\infty} \frac{x^k(F_z x^k + F_{z+1})}{k\,(1 - x^k - x^{2k})} = \frac{x\,U'(x)}{U(x)} + \sum_{k=1}^{\infty} G(k,x)$$

where

$$G(k,x) = \frac{x^k(F_z x^k(-(5k+1)x^{2k} + (k+1)x^{3k} + 3(k-1)x^k - 4k + 2) - F_{z+1}(6kx^{2k} - (k-1)x^{3k} + (k+1)x^{4k} + (k+1)x^k + k - 1))}{(x^{2k} + x^k - 1)^3}$$

$$b(r_n) = n + G(1, r_n) + \sum_{k=2}^{\infty} G(k, r_n)$$

Now for $x = r_n$ is (independently on $z$)

$$\lim_{k \to \infty} G(k, r_n)^{1/k} = \frac{1}{\varphi} < 1$$

and the sum tends to a constant as $n$ tends to infinity.

$$\sum_{k=2}^{\infty} G(k, r_n) = c$$

For example if $z = 0$ then $c = 19.55999649742693171136312985683987556 3\ldots$

If $k = 1$ then we obtain

$$G(1, r_n) = \frac{2\,r_n^2\left(((r_n - 3)\,r_n^2 - 1)F_z - (r_n^3 + 3r_n + 1)F_{z+1}\right)}{(r_n^2 + r_n - 1)^3} \sim 2\,\varphi^{-z/2}\,5^{1/4}n^{3/2}$$

Together

$$b(r_n) \sim n + 2\,\varphi^{-z/2}\,5^{1/4}n^{3/2} + c \sim 2\,\varphi^{-z/2}\,5^{1/4}n^{3/2}$$

$$U(r_n) = e^{\log(U(r_n))} = e^{\sum_{k=1}^{\infty} \frac{r_n^k(F_z r_n^k + F_{z+1})}{k(1 - r_n^k - r_n^{2k})}}$$

We have

$$\sum_{k=1}^{\infty} \frac{r_n^k(F_z r_n^k + F_{z+1})}{k(1 - r_n^k - r_n^{2k})} = \frac{r_n(F_z r_n + F_{z+1})}{1 - r_n - r_n^2} + \sum_{k=2}^{\infty} \frac{r_n^k(F_z r_n^k + F_{z+1})}{k(1 - r_n^k - r_n^{2k})}$$

Contribution of the first term is

$$\frac{r_n(F_z r_n + F_{z+1})}{1 - r_n - r_n^2} \sim -\frac{F_{z+1}}{10} + \frac{1}{20}(\sqrt{5} - 9)F_z + \frac{\varphi^{z/2}}{5^{1/4}}\sqrt{n} = \left(\frac{\phi}{10} - \frac{1}{2}\right)F_z - \frac{F_{z+1}}{10} + \frac{\varphi^{z/2}}{5^{1/4}}\sqrt{n}$$

For $k > 1$

$$\frac{r_n^k(F_z r_n^k + F_{z+1})}{k(1 - r_n^k - r_n^{2k})} \sim -\frac{2^k\left((\sqrt{5}+1)^k F_{z+1} + 2^k F_z\right)}{\left((2(\sqrt{5}+1))^k + 4^k - (\sqrt{5}+1)^{2k}\right)k} = \frac{F_{z+1}\varphi^k + F_z}{k(\varphi^{2k} - \varphi^k - 1)}$$

$$U(r_n) = e^{\sum_{k=1}^{\infty} \frac{r_n^k(F_z r_n^k + F_{z+1})}{k(1 - r_n^k - r_n^{2k})}} \sim e^{\left(\frac{\phi}{10} - \frac{1}{2}\right)F_z - \frac{F_{z+1}}{10} + \frac{\varphi^{z/2}}{5^{1/4}}\sqrt{n} + \sum_{k=2}^{\infty} \frac{F_{z+1}\varphi^k + F_z}{k(\varphi^{2k} - \varphi^k - 1)}}$$

Together

$$a_n \sim \frac{U(r_n)}{\sqrt{2\pi * b(r_n)}} \frac{1}{r_n^n} = \frac{e^{\left(\frac{\phi}{10} - \frac{1}{2}\right)F_z - \frac{F_{z+1}}{10} + \frac{\varphi^{z/2}}{5^{1/4}}\sqrt{n} + \sum_{k=2}^{\infty} \frac{F_{z+1}\varphi^k + F_z}{k(\varphi^{2k} - \varphi^k - 1)}}}{\sqrt{2\pi * (2\,\varphi^{-z/2}\,5^{1/4} n^{3/2})}} * \varphi^n \exp\left(\frac{\varphi^{z/2}\sqrt{n}}{5^{1/4}}\right)$$

The final asymptotic is

$$a_n \sim \frac{\varphi^{n + \frac{z}{4}} \exp\left(\left(\frac{\varphi}{10} - \frac{1}{2}\right)F_z - \frac{1}{10}F_{z+1} + \frac{2\,\varphi^{z/2}}{5^{1/4}}\sqrt{n} + \sum_{k=2}^{\infty} \frac{F_z + F_{z+1}\varphi^k}{(\varphi^{2k} - \varphi^k - 1)k}\right)}{2\sqrt{\pi}\,5^{1/8}\,n^{3/4}}$$

# Applications

The sequence A166861 ($z = 0$)
Generating function:

$$\prod_{k=1}^{\infty} \frac{1}{(1-x^k)^{F_k}}$$

Asymptotics:

$$a_n \sim \frac{\varphi^n \; e^{-\frac{1}{10} + \frac{2}{5^{1/4}} \sqrt{n} + S}}{2\sqrt{\pi} \; 5^{1/8} \; n^{3/4}}$$

where

$$S = \sum_{k=2}^{\infty} \frac{\varphi^k}{(\varphi^{2k} - \varphi^k - 1)\,k} = \sum_{k=2}^{\infty} \frac{1}{2k\sinh(k\,\operatorname{arccsch}(2)) - k}$$

$S = 0.60047660139257591296971949485039357608376512393964351135554 7131467\ldots$

Numerical verification (20000 terms), the ratio tends to 1:

```
z = 0; s = NSum[(Fibonacci[z] + Fibonacci[z+1] * GoldenRatio^k) / ((GoldenRatio^(2*k) - 
GoldenRatio^k - 1)*k), {k, 2, Infinity}, WorkingPrecision -> 100, AccuracyGoal -> 100, 
PrecisionGoal -> 100, NSumTerms -> 10000];
Show[Plot[1, {n, 1, 20000}, PlotStyle -> Red], 
ListPlot[Table[A166861[[n]]/((GoldenRatio^(n + z/4) / (2 * Sqrt[Pi] * 5^(1/8) * 
n^(3/4))) * Exp[(GoldenRatio/10 - 1/2)*Fibonacci[z] - Fibonacci[z + 1]/10 + 
(2*GoldenRatio^(z/2) * Sqrt[n])/5^(1/4) + s]), {n, 1, Length[A166861]}]], PlotRange -> 
{0.5, 1}, AxesOrigin -> {0, 0.5}]
```

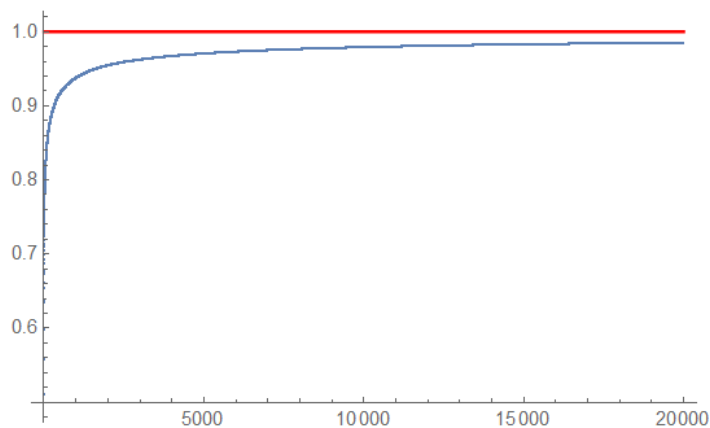

The sequence A200544 ($z = 1$)

Generating function:
$$\prod_{k=1}^{\infty} \frac{1}{(1-x^k)^{F_{k+1}}}$$

Asymptotics:
$$a_n \sim \frac{\varphi^{n+\frac{1}{4}} \, e^{\frac{\varphi}{10} - \frac{3}{5} + \frac{2\varphi^{1/2}}{5^{1/4}} \sqrt{n} + S}}{2\sqrt{\pi} \, 5^{1/8} \, n^{3/4}}$$

where
$$S = \sum_{k=2}^{\infty} \frac{1+\varphi^k}{(\varphi^{2k} - \varphi^k - 1)\,k}$$

$S = 0.7902214013751085262994702391769374769675268259229550490716908010345\ldots$

---

The sequence A260787 ($z = 2$)

Generating function:
$$\prod_{k=1}^{\infty} \frac{1}{(1-x^k)^{F_{k+2}}}$$

Asymptotics:
$$a_n \sim \frac{\varphi^{n+\frac{1}{2}} \, e^{\frac{\varphi}{10} - \frac{7}{10} + \frac{2\varphi}{5^{1/4}} \sqrt{n} + S}}{2\sqrt{\pi} \, 5^{1/8} \, n^{3/4}}$$

where
$$S = \sum_{k=2}^{\infty} \frac{1 + 2\varphi^k}{(\varphi^{2k} - \varphi^k - 1)\,k}$$

$S = 1.390698002767684439269189734027331053051291949862598560427237 9325\ldots$

If we set $F_{-1} = 1$ (according to Mathematica), then my formula is correct also for $z = -1$, the sequence A109509. Generating function:

$$\prod_{k=1}^{\infty} \frac{1}{(1-x^k)^{F_{k-1}}}$$

Asymptotics:

$$a_n \sim \frac{\varphi^{n-\frac{1}{4}} e^{\frac{\varphi}{10} - \frac{1}{2} + \frac{2\varphi^{-1/2}}{5^{1/4}}\sqrt{n} + S}}{2\sqrt{\pi}\, 5^{1/8}\, n^{3/4}}$$

where

$$S = \sum_{k=2}^{\infty} \frac{1}{(\varphi^{2k} - \varphi^k - 1)\, k}$$

$S =$ 0.189744799982532613329750744326543900883761701983311537716143669 ...

---

## References:

[1] OEIS - The On-Line Encyclopedia of Integer Sequences

[2] A. M. Odlyzko, Asymptotic enumeration methods, pp. 1063-1229 of R. L. Graham et al., eds., Handbook of Combinatorics, 1995

[3] V. Kotěšovec, Asymptotics of sequence A034691, Sep 09 2014, mirror

[4] The website http://web.telecom.cz/vaclav.kotesovec/math.htm

---